\documentclass[11pt,a4paper,twocolumn]{article}
\usepackage[utf8]{inputenc}
\usepackage{amsfonts}
\usepackage{graphicx}      
\usepackage[square,sort,comma,numbers]{natbib}
\usepackage{color}
\usepackage{amsmath, amsthm, amsfonts, amssymb, psfrag}
\usepackage{authblk}
\usepackage[left=1.50cm, right=1.50cm, top=2.00cm, bottom=2.00cm]{geometry}

\usepackage{sectsty}
\sectionfont{\fontsize{11}{11}\selectfont}
\subsectionfont{\fontsize{11}{11}\selectfont}

\begin{document}
	
	\title{A Multiobjective MPC Approach for Autonomously Driven Electric Vehicles}
	\author[*]{Sebastian Peitz}
	\author[*]{Kai Sch{\"a}fer}
	\author[**]{Sina Ober-Bl{\"o}baum}
	\author[***]{Julian Eckstein}
	\author[***]{Ulrich K{\"o}hler}
	\author[*]{Michael Dellnitz}
	\affil[*]{\normalsize Department of Mathematics, Paderborn University, Warburger Str.~100, 33098 Paderborn, Germany}
	\affil[**]{Department of Engineering Science, University of Oxford, Parks Road, Oxford OXI 3PJ, UK}
	\affil[***]{Hella KGaA Hueck \& Co., Beckumer Str. 130, 59552 Lippstadt, Germany}

	\date{\vspace{-5ex}}

	\twocolumn[
	\begin{@twocolumnfalse}	
	
	\maketitle
	
	\begin{abstract}                
		We present a new algorithm for model predictive control of non-linear systems with respect to multiple, conflicting objectives. The idea is to provide a possibility to change the objective in real-time, e.g.~as a reaction to changes in the environment or the system state itself. The algorithm utilises elements from various well-established concepts, namely multiobjective optimal control, economic as well as explicit model predictive control and motion planning with motion primitives. In order to realise real-time applicability, we split the computation into an online and an offline phase and we utilise symmetries in the open-loop optimal control problem to reduce the number of multiobjective optimal control problems that need to be solved in the offline phase. The results are illustrated using the example of an electric vehicle where the longitudinal dynamics are controlled with respect to the concurrent objectives arrival time and energy consumption.
	\end{abstract}
	\vspace{2ex}
	
	\end{@twocolumnfalse}
	]
	
	
	\section{Introduction}
	\label{sec:Introduction}
	In many applications from industry and economy, the simultaneous optimisation of several criteria is of great interest. In transportation, for example, one wants to reach a destination as fast as possible while minimising the energy consumption. This example illustrates that in general, the different objectives contradict each other. Therefore, the task of computing the set of optimal compromises between the conflicting objectives, the so-called \emph{Pareto set}, arises, leading to a multiobjective optimisation problem (MOP) or multiobjective optimal control problem (MOCP). Based on the knowledge of the Pareto set, a \emph{decision maker} can design improved systems or even allow for changes in control parameters during operation as a reaction on external influences or changes in the system state itself. 
	There exist various algorithms for the solution of MOCPs such as scalarisation techniques (cf.~\cite{Ehr05} for an overview), evolutionary algorithms (\cite{CLV07}) or set oriented methods (\cite{SWO+13}). All approaches have in common that a large number of function evaluations is typically needed. Thus, the direct computation of the Pareto set is time consuming and a computation in real-time is not possible. 
	However, in particular the design of optimal drive strategies requires online adaption of control strategies. This is even more the case now that autonomous driving and battery electric vehicles (EVs) with comparatively low ranges are both gaining increased attention, requiring advanced control algorithms.
	
	Control theory has been influenced significantly by the advances in computational power during the last decades. For a large variety of systems, it is nowadays possible to use model based optimal control algorithms to design sophisticated feedback laws. This concept is known as model predictive control (MPC) (see e.g.~\cite{Mac02, GP11}). The general goal of MPC is to stabilise a system by using a combination of open and closed-loop control: using a model of the system dynamics, an open-loop optimal control problem is solved in real-time over a so-called \emph{prediction horizon}. The first part of this solution is then applied to the real system while the optimisation is repeated to find a new control function, with the prediction horizon moving forward (for this reason, MPC is also referred to as moving horizon control or receding horizon control).
	
	Due to the huge success of MPC, a large variety of algorithms has been established, where a first distinction can be made between linear and non-linear MPC. The first category refers to schemes in which linear models and quadratic objective functions are used to predict the system dynamics. The resulting optimisation problems are convex, i.e.~global solutions can be computed very fast. Linear MPC approaches have been very successful in a large variety of industrial applications (see e.g.~\cite{QB97} and \cite{LC97} for an overview in applications and theory). The advantage of non-linear MPC (\cite{GP11}), on the other hand, is that the typically non-linear system behaviour can be approximated in a more accurate way. Furthermore, special optimality criteria and non-linear constraints can be incorporated easily. However, the complexity and thus the time to solve the resulting optimisation problem increases such that it is often difficult to preserve real-time capability (see e.g.~\cite{EPS+16}). Further extensions are, for example, \emph{economic MPC} (see e.g.~\cite{RA09, DAR11}) or \emph{explicit MPC} (see e.g.~\cite{AB09}). In the first approach alternative, \emph{economic} objectives are pursued instead of stabilising the system. In the second approach the problem of real-time applicability is addressed by introducing an offline phase during which the open-loop optimal control problem is solved for a large number of possible situations, using e.g.~multi-parametric non-linear programming. The solutions are then stored in a library such that they are directly available in the online phase.
	
	Another way for optimal strategy planning is the concept \emph{motion planning with motion primitives} going back to \cite{FDF05} (see also \cite{Kob08, FOK12}). The challenge of online applicability is addressed with a two-phase approach similar to explicit MPC but here, valid control as well as state trajectories are obtained by combining several short pieces of simply controlled trajectories that are stored in a motion planning library. These motion primitives can be sequenced to longer trajectories in various combinations. In the online phase, the optimal sequence of motion primitives is determined from the motion planning library using e.g.~graph search methods (see e.g~\cite{Kob08}). To reduce the computational effort, the motion primitive approach extensively relies on exploiting symmetries in the dynamical control system such that a motion primitive can be used in multiple situations, e.g.~by performing a translation or rotation under which the dynamics are invariant.
	
	In this article, we present a new algorithm for multiobjective MPC of non-linear systems. Problems with multiple criteria have been addressed by several authors using scalarisation techniques (see e.g.~\cite{BP09} for a weighted sum or \cite{ZFT12} for a reference point approach). For non-convex problems, scalarisation approaches often face difficulties such that we here want to compute the entire Pareto set in advance. To this end, we combine elements from multiobjective optimal control, explicit MPC and motion planning with motion primitives. The resulting algorithm consists of an offline phase during which multiobjective optimal control problems are solved and stored in a library for a wide range of possible scenarios (i.e.~constant velocity, braking, accelerating). Invariances in the optimal control problem are exploited in order to reduce the number of problems that need to be solved. In the online phase, the currently active scenario is identified and the corresponding Pareto set is selected from the library. According to a decision maker's preference, an optimal compromise is then selected from the Pareto set and the first part of the solution is applied to the system. Similar to MPC, this is done repeatedly such that a feedback control behaviour is realised. The difference to other approaches is the possibility to interactively choose between different objectives such that the system behaviour can be modified easily. This can be very useful for autonomous driving, where one is interested in reaching a destination as fast as possible while minimising the energy consumption.
	
	The outline of the article is as follows. In Section~\ref{sec:Problem_formulation}, we introduce the multiobjective MPC problem and the concept of Pareto optimality before describing the algorithm in detail and comparing it to other MPC approaches. In Section~\ref{sec:Application_EV}, we describe the application of the algorithm to an electric vehicle. The aim is to realise autonomous driving where the passenger can decide between the objectives fast and energy efficient driving. We present the results in Section~\ref{sec:Results} before drawing a conclusion in Section~\ref{sec:Conclusion}.
	
	\section{Problem Formulation and Methodology}
	\label{sec:Problem_formulation}
	Before describing the algorithm, we will briefly introduce the two main concepts we will be making use of, namely multiobjective optimal control and model predictive control. For more detailed introductions, we refer to \cite{Ehr05} and \cite{GP11}, respectively.
	
	A \emph{multiobjective optimal control} problem (MOCP) can be formulated mathematically using differential(-algebraic) equations describing the physical behaviour of the system together with optimisation criteria and optimisation constraints in the following way
	\begin{equation}
	\min_{x,u,t_f} J(x,u,t_f) = \int_{t_0}^{t_f} C(x(t),u(t))\, dt + \Phi(x(t_f))\label{eq:J}
	\end{equation}
	such that
	\begin{align}
	& \dot{x}(t)  = f(x(t),u(t))\quad \forall t \in [t_0,t_f],\quad x(t_0)=x_0 \label{eq:diff}\\
	& h(x(t),u(t)) \le 0 \quad \forall t \in [t_0,t_f],\label{eq:constraints}
	\end{align}
	where $x(t)\in \mathcal{X}$ is the system state (e.g. the position and velocity of a car) and $u(t)\in \mathcal{U}$ the control (e.g.~the engine torque or the steering wheel position). $\mathcal{X}$ and $\mathcal{U}$ are the spaces of feasible states and controls, respectively. 
	The constraints may depend on the state as well as the control, e.g.~limiting the velocity or energy consumption. $J$ describes criteria that have to be optimised. When there exists a unique solution $x(t) \in \mathcal{X}$ for every $u(t) \in \mathcal{U}$ and $x_0 \in \mathcal{X}$ and we fix the time frame, we can introduce a reduced objective $J: \mathcal{U} \times \mathcal{X} \rightarrow \mathbb{R}^k$, where $k$ is the number of objectives, and the corresponding reduced problem:
	\begin{equation}
	\min_{u} J(u, x_0) = \int_{t_0}^{t_f} C(\varphi_u(x_0,t))\, dt + \Phi(\varphi_u(x_0,t_f)). \label{eq:Jred}
	\end{equation}
	Here $\varphi_u(x_0,t)$ is the flow of the dynamical control system \eqref{eq:diff}.
	
	\begin{figure}[h!]
		\centering
		\parbox[b]{0.24\textwidth}{\centering \includegraphics[width=0.2\textwidth]{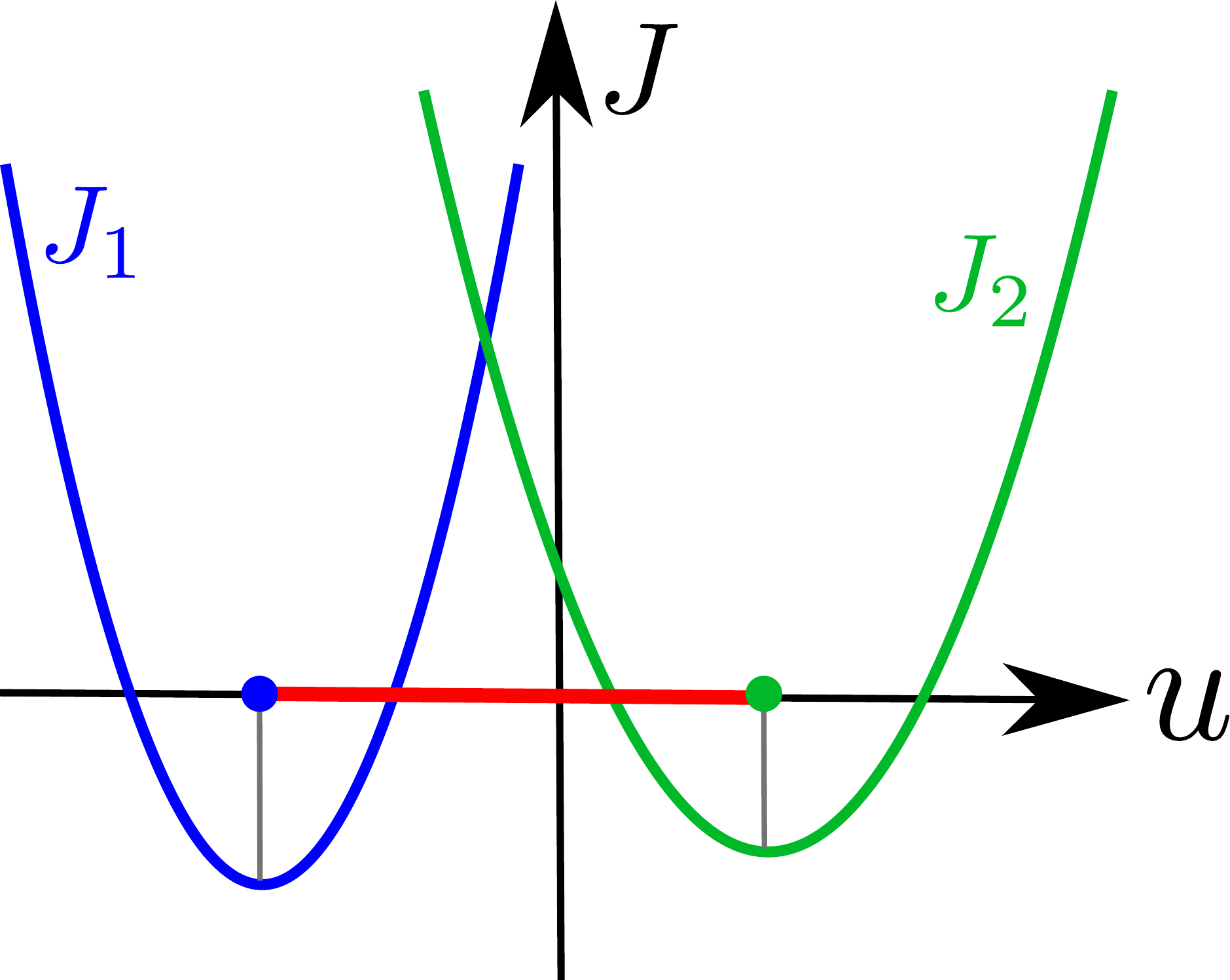}\\(a)}
		\parbox[b]{0.24\textwidth}{\centering \includegraphics[width=0.18\textwidth]{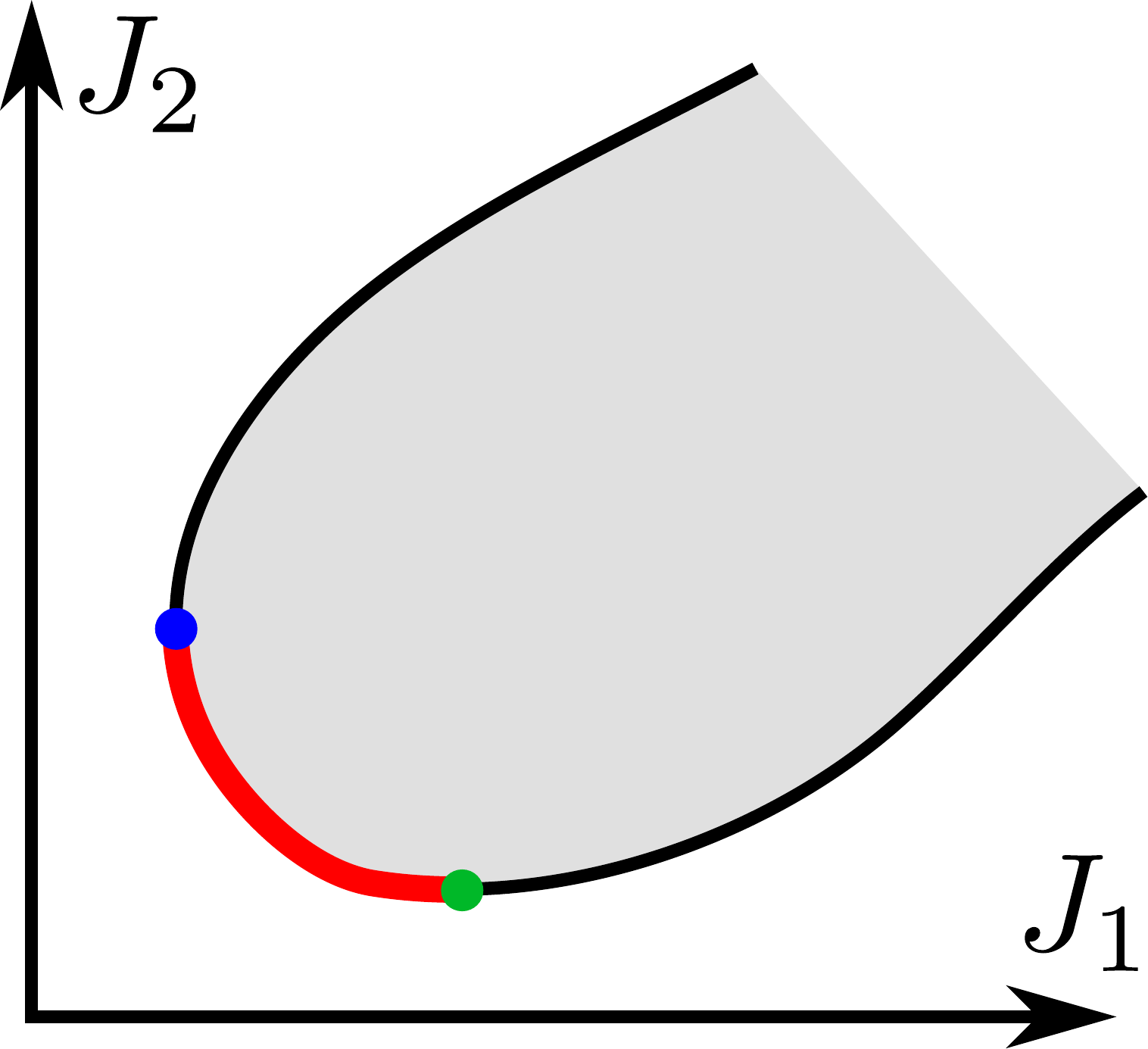}\\(b)}
		\caption{Pareto set (a) and front (b) of the multiobjective optimisation problem $\min_{u\in\mathbb{R}} J(u)$, $J: \mathbb{R}\rightarrow\mathbb{R}^2$.}
		\label{fig:MOP}
	\end{figure}
	In many applications from industry and economy, one is interested in simultaneously optimising not only one but \emph{several} criteria and hence, $k>1$ and $J$ is vector-valued. In this situation the solution does in general not consist of isolated optimal points but of the \emph{set of optimal compromises}, the so-called \emph{Pareto set} (cf.~\cite{Ehr05} for a detailed introduction). The set consists of all functions $u(t)$ that are \emph{nondominated}, i.e.~for which there does not exist a solution $u^*(t)$ that is superior in all objectives (cf.~Figure~\ref{fig:MOP}).
	\begin{figure}[h!]
		\centering
		\parbox[b]{0.24\textwidth}{\centering \includegraphics[width=0.24\textwidth]{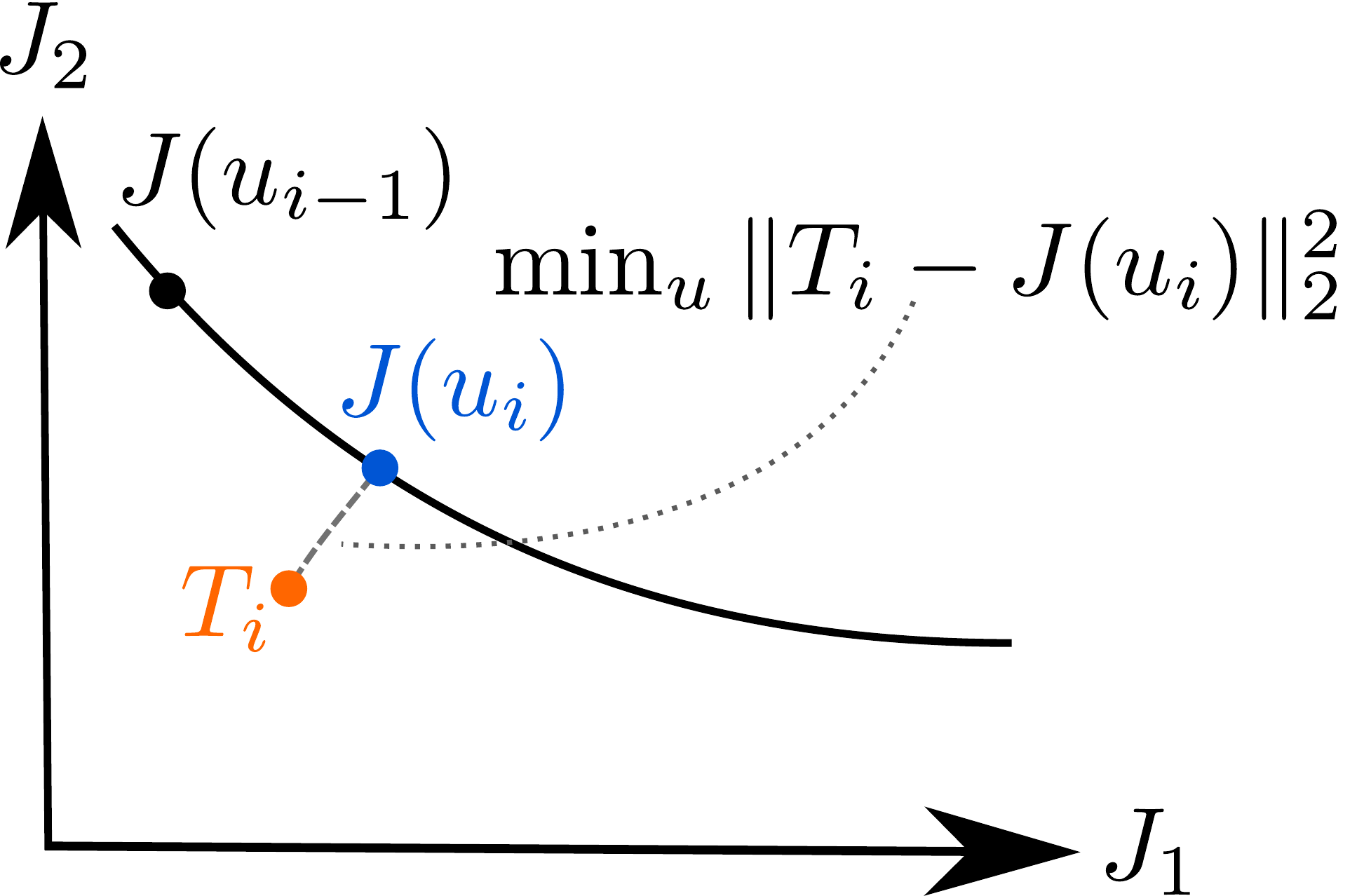}\\(a)}
		\parbox[b]{0.24\textwidth}{\centering \includegraphics[width=0.21\textwidth]{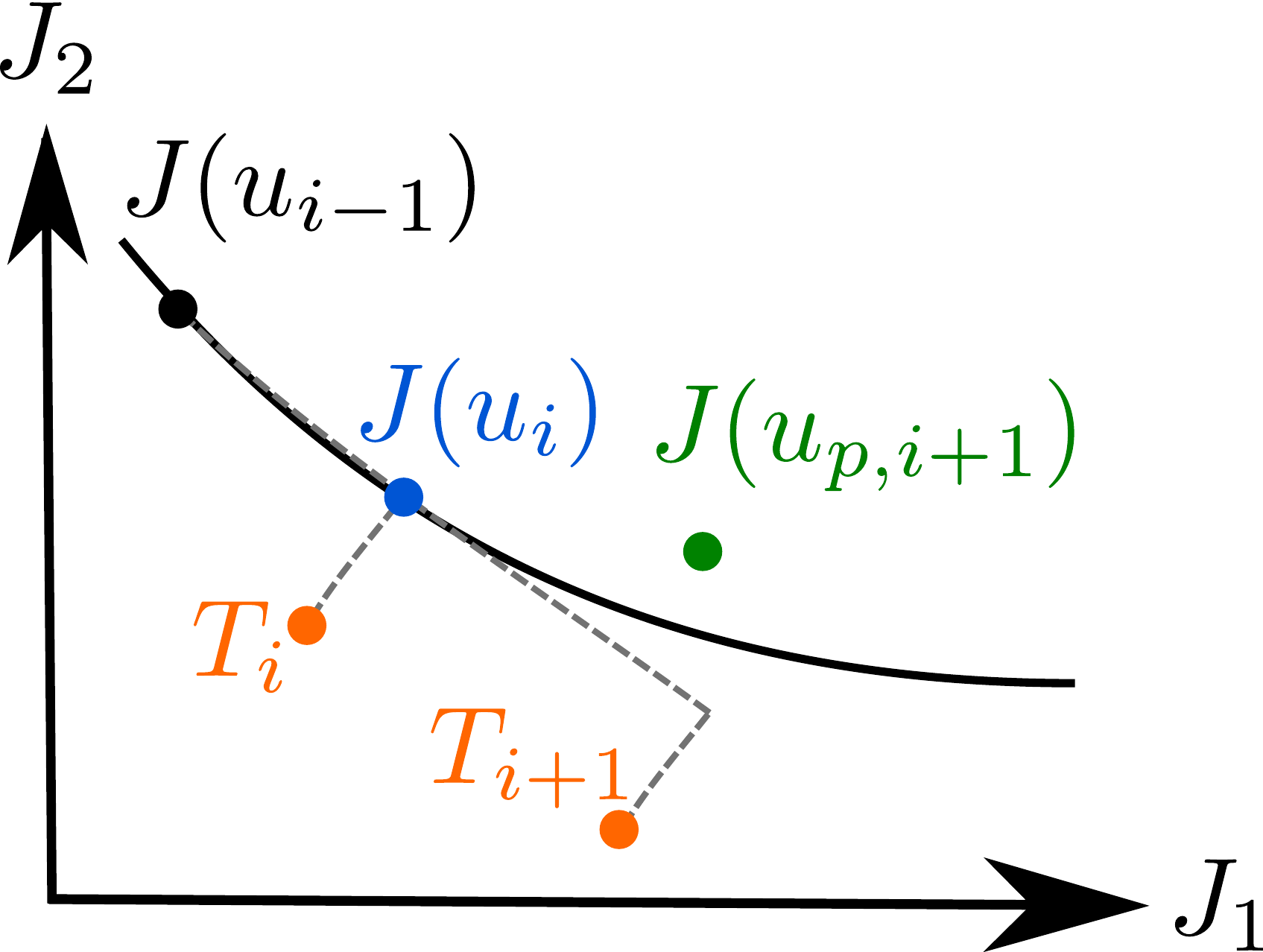}\\(b)}
		\caption{Reference point method in image space. (a) Determination of the $i$-th point on the Pareto front by solving a scalar optimisation problem. (b) Computation of new target point $T_{i+1}$ and predictor step in decision space ($u_{p,i+1})$.}
		\label{fig:ReferencePoint}
	\end{figure}
	For the solution of \eqref{eq:Jred}, we here use a scalarisation technique by which the Pareto set is approximated by a finite set of points that are computed consecutively by minimising the euclidean distance between a point $J(u, x_0)$ and a so-called \emph{target point} $T$ which lies outside the reachable set in image space (see Figure~\ref{fig:ReferencePoint} for an illustration). Since a point computed this way lies on the boundary of the reachable set, there exists no point which is superior in all objectives and hence, the point is Pareto optimal. Starting with one point (e.g.~the scalar minimum of one of the objectives), the next points can be computed recursively until the other end of the Pareto front (i.e. the other scalar minimum) is reached. In \cite{DEF+16}, this method is used to compute the Pareto set for the conflicting objectives driven distance and energy consumption for EVs. The scalar optimal control problems are solved using an SQP method (cf.~\cite{NW06}). 
	
	\begin{figure}[h!]
		\centering
		\includegraphics[width=0.4\textwidth]{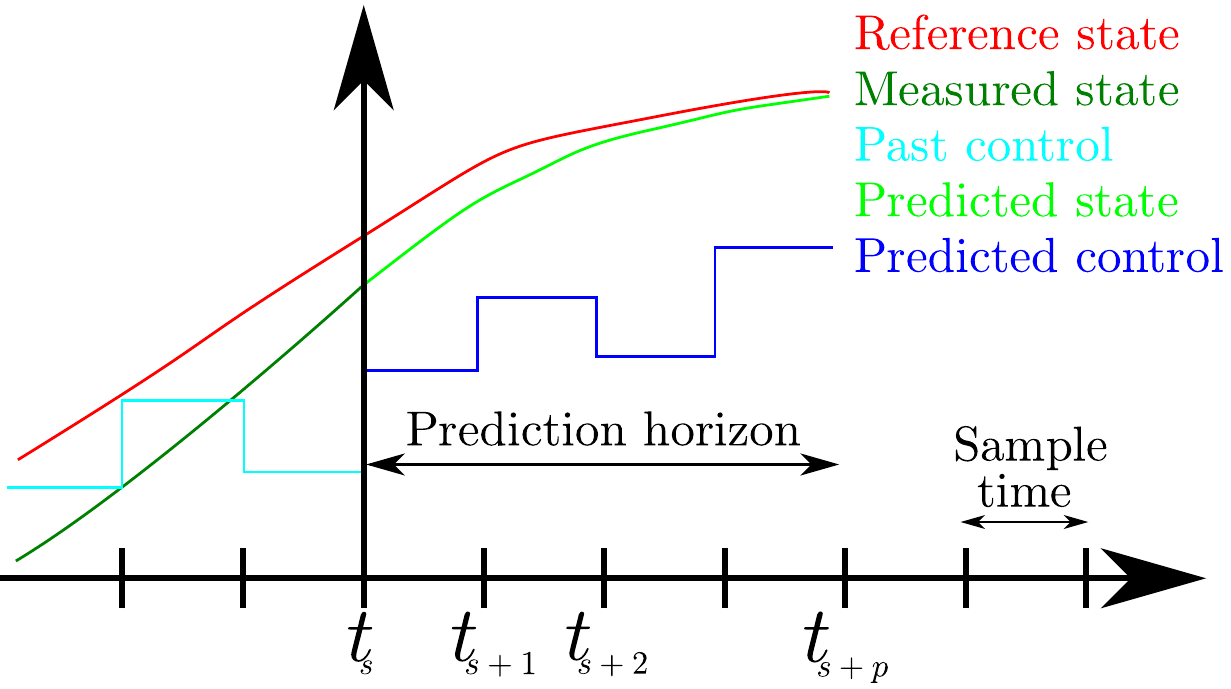}
		\caption{Sketch of the MPC methodology. While the first part of the predicted control is applied to the system, the next control is predicted (via open-loop optimal control) on a shifted horizon.}
		\label{fig:MPC}
	\end{figure}
	The algorithm presented here builds on these results, but we need to extend them in order to construct a feedback controller. This is realised by an \emph{MPC} approach, where the problem \eqref{eq:Jred} is solved repeatedly for varying time frames ($t_0 = t_s$, $t_f = t_{s+p}$, $s = 1,2,\ldots$) while the system is running at the same time. Then, the first interval of the predicted control, $u(t_s)$, is applied to the real system and the optimal control problem is solved again with a time frame shifted by one. The procedure is illustrated in Figure~\ref{fig:MPC}. The concept of MPC was initially developed to stabilise a system (\cite{GP11}), i.e.~to drive the system state to a (potentially time dependent) reference state. However, stabilisation is not always the main concern. Considering the EV, for example, we only require a part of the state, namely the velocity, to remain within prescribed bounds, which then gives us the opportunity to pursue additional objectives such as minimising the energy consumption. This concept is known as \emph{economic MPC} (see e.g.~\cite{RA09,DAR11}).
	
	\subsection{The Offline-Online Multiobjective MPC Concept}
	\label{subsec:Algorithm}
	Since MOCPs are considerably more expensive to solve than scalar problems, it is computationally infeasible to directly include them in an MPC framework. A simple way to circumvent this problem is to scalarise the objective function by introducing a weighting factor (i.e.~$\widehat{J}=\sum_{i=1}^k \rho_i J_i, \rho_i\in[0,1]$). In this case however, an assumption has to be made in advance which can in practice lead to unfavourable results. A slight increase in one objective might allow for a strong reduction in another one, for example. Hence, we are interested in providing the entire Pareto set during the MPC routine. To avoid large computing times during execution, we therefore split the computation in an \emph{offline} and an \emph{online phase}, similar to explicit MPC approaches (cf.~\cite{AB09}). 
	
	The \emph{offline phase} consists of several steps. First, various \emph{scenarios} are identified for which MOCPs need to be solved. The scenarios are determined by the system states and the constraints. Secondly, in order to reduce the number of scenarios, the dynamical control system is analysed with respect to invariances, which are formally described by a finite-dimensional Lie group $G$ and its group action $\psi: \mathcal{X} \times G \rightarrow \mathcal{X}$. 
	A dynamical control system, described by \eqref{eq:diff}, is invariant under the group action $\psi$, or equivalently, $G$ is a symmetry group for the system \eqref{eq:diff}, if for all $g\in G$, $x_0\in \mathcal{X}$, $t\in [t_0,t_f]$ and all piecewise-continuous control functions $u:[t_0,t_f]\rightarrow \mathcal{U}$ it holds
	\begin{equation}
	\psi(g,\varphi_u(x_0,t)) = \varphi_u(\psi(g,x_0),t) \quad \forall g\in G. \label{eq:Invariance}
	\end{equation}
	That means that the group action on the state commutes with the flow. 
	Invariance leads to the concept of {\em equivalent trajectories}. 
	Two trajectories are equivalent if they can be exactly superimposed through time translation and the action of the symmetry group. 
	In the classical concept of motion primitives (\cite{FDF05}), all equivalent trajectories are summed up in an equivalence class, i.e.~only a single representative is stored that can be used at many different points when transformed by the symmetry action. 
	In other words, controlled trajectories that have been computed for a specific situation are suitable in many different (equivalent) situations as well.
	In our approach, we extend this concept by identifying symmetries in the solution of the MOCP \eqref{eq:Jred} with respect to the initial conditions $x_0$:
	\begin{align}
	\arg \min_{u} J(u, x_0) = \arg \min_{u} J(u, \psi(g,x_0)) \quad \forall g \in G. \label{eq:Invariance_MOCP}
	\end{align}
	This means that we require the \emph{Pareto set} to be invariant under group actions on the initial conditions. If the objective function is also invariant under the same group action, then all trajectories contained in an equivalence class defined by \eqref{eq:Invariance} will also be contained in an equivalence class defined by \eqref{eq:Invariance_MOCP}. However, this class may contain more solutions since we do not explicitly pose restrictions on the state but only require the solution of \eqref{eq:Jred} to be identical. Alternatively, if the objective function is linear in the states and the group action corresponds to translations in initial states, we do not require invariance of the objective function to satisfy \eqref{eq:Invariance_MOCP}.
	
	Identifying invariances according to \eqref{eq:Invariance_MOCP}, the number of MOCPs can be reduced. If the system is invariant under translation of the initial position $p(t_0)$, for example, we do not need to solve multiple MOCPs that only differ in the position. Once these equivalence classes have been identified, we can reduce the number of possible scenarios accordingly. We then solve the resulting MOCPs on the prediction horizon $T_p$, introduce a parametrisation $\rho$ (which can then be chosen by the decision maker in the online phase) and store the Pareto sets and fronts in a library such that they can be used in the online phase. Since in general there is an infinite number of feasible initial conditions, there consequently exists an infinite number of scenarios that we have to consider. In practice, this obviously cannot be realised and we have to introduce a finite set of scenarios. In the online phase, we then pick the scenario that is closest to the true initial condition. If a violation of the state constraints has to be avoided (the EV, e.g., is not allowed to go faster than the maximum speed), then a selection towards the ''safe'' side can be made. In case of the EV, we would consequently pick a solution corresponding to a velocity slightly higher than the actual velocity. This way, the maximally allowed acceleration would be bounded such that exceeding the speed limit is not possible.
	
	The \emph{online phase} is now basically a standard MPC approach, the difference being that we obtain the solution of our control problem from a library instead of solving it in real-time, similar to explicit MPC approaches:
	\begin{itemize}
		\item[1.] measure the current system states that are necessary for the identification of the current scenario,
		\item[2.] choose the corresponding Pareto set from the library, i.e.~the one with initial conditions closest to the current system state. (Due to the approximation, we cannot formally guarantee that the constraints are not violated. However, as a start we consider applications where this is acceptable.)
		\item[3.] choose one optimal compromise $u$ from the set, according to a decision maker's preference $\rho$,
		\item[4.] apply the first step (i.e.~the sample time) of the solution $u$ to the real system and go back to 1.
	\end{itemize}
	The resulting algorithm thus provides a feedback law. In the offline phase, we define the scenarios in such a manner that the system cannot be steered out of the set of feasible states. This means that only controls $u$ are valid that do not lead to a violation of the constraints. Additionally, we include scenarios which steer the system into the set of feasible states from any initial condition. In the literature, this is known as \emph{viability}, cf.~\cite{GP11}. In case of the EV, for example, we have to include controls such that the velocity can be steered to values satisfying the constrains from any initial velocity.
	
	The presented algorithm can be seen as an extension of (extended) MPC approaches to multiple objectives. We consider \emph{economic} objectives (cf.~\cite{RA09}) and do not focus on the stabilisation of the system. This allows us to pursue multiple objectives between which a decision maker can choose dynamically, e.g.~in order to react on changes in the environment or the system state itself. In contrast to weighting methods, the entire Pareto set is known, providing increased system knowledge.
	
	\section{Application to Electric Vehicle}
	\label{sec:Application_EV}
	In this section the algorithm is utilised to control the longitudinal dynamics of an EV, thereby extending prior work, see \cite{DEF+14} for a scalar optimal control problem, \cite{DEF+16} for a multiobjective optimal control problem and \cite{EPS+16} for a comparison of two scalar MPC approaches.
	
	\subsection{Vehicle Model}
	\label{subsec:EV_model}
	The EV model is derived by coupling the equations for the electrical and the mechanical subsystem via efficiency maps. This yields a system of four coupled, non-linear ordinary differential equations for the system state $x(t) = \left( v(t), S(t), U_{d,L}(t), U_{d,S}(t) \right)$. Here, $v$ is the vehicle velocity, $S$ is the battery state of charge and $U_{d,L}$ and $U_{d,S}$ are the long and short term voltage drops, respectively. The system is controlled by setting the torque $u(t)$ of the front wheels. Additionally, the battery current $I(t)$ is computed from the state $x(t)$ via an algebraic equation and the position by integrating the velocity: $p(t) = \int_{t_0}^{t}v(\tau)d\tau$. For the derivation and the exact formulation of the dynamical system, we refer the reader to \cite{EPS+16}.
	
	Based on the system dynamics, we formulate the MOCP for the EV with variable final time:
	\begin{align}
	&min_{u} \left( \begin{array}{c} \label{eq:MOCP_EV_J}
	S(t_0)-S(t_f) \\ t_f - t_0
	\end{array} \right), \\
	\dot{x}(t)&=f(x(t),u(t)), \label{eq:MOCP_EV_diff}\\
	v_{min}(t) &\le v(t) \le v_{max}(t), \quad &t \in [0,t_f] \label{eq:MOCP_EV_constraints1} \\
	I_{min}(t) &\le I(t) \le I_{max}(t), \quad &t \in [0,t_f] \label{eq:MOCP_EV_constraints2} \\
	x(0) &= x_0, \ p(t_f) = p_f. \label{eq:MOCP_EV_x0}
	\end{align}
	We set the final position $p_f$ to $100$\,m, which means that we here define the prediction horizon based on the position. Correspondingly, the sample time is also specified with respect to the position, $\delta = 20$\,m. The conflicting objectives are to reach $p_f$ as fast as possible ($J_2$) while minimising the energy consumption ($J_1$). The battery current $I$ is limited in order to avoid damaging the battery which results in implicit constraints on the control $u$. The velocity constraints are part of the scenarios which are defined in the offline phase.
	
	\subsection{Offline Phase: System Analysis and Solution of Multi- objective Optimal Control Problems}
	\label{subsec:EV_offline_phase}
	In this section we describe how the different steps of the offline phase are applied to the EV.
	\subsubsection{Symmetry Analysis}
	\label{subsec:EV_symmetries}
	\begin{figure}[h!]
		\centering
		\parbox[b]{0.4\textwidth}{\centering \includegraphics[width=0.4\textwidth]{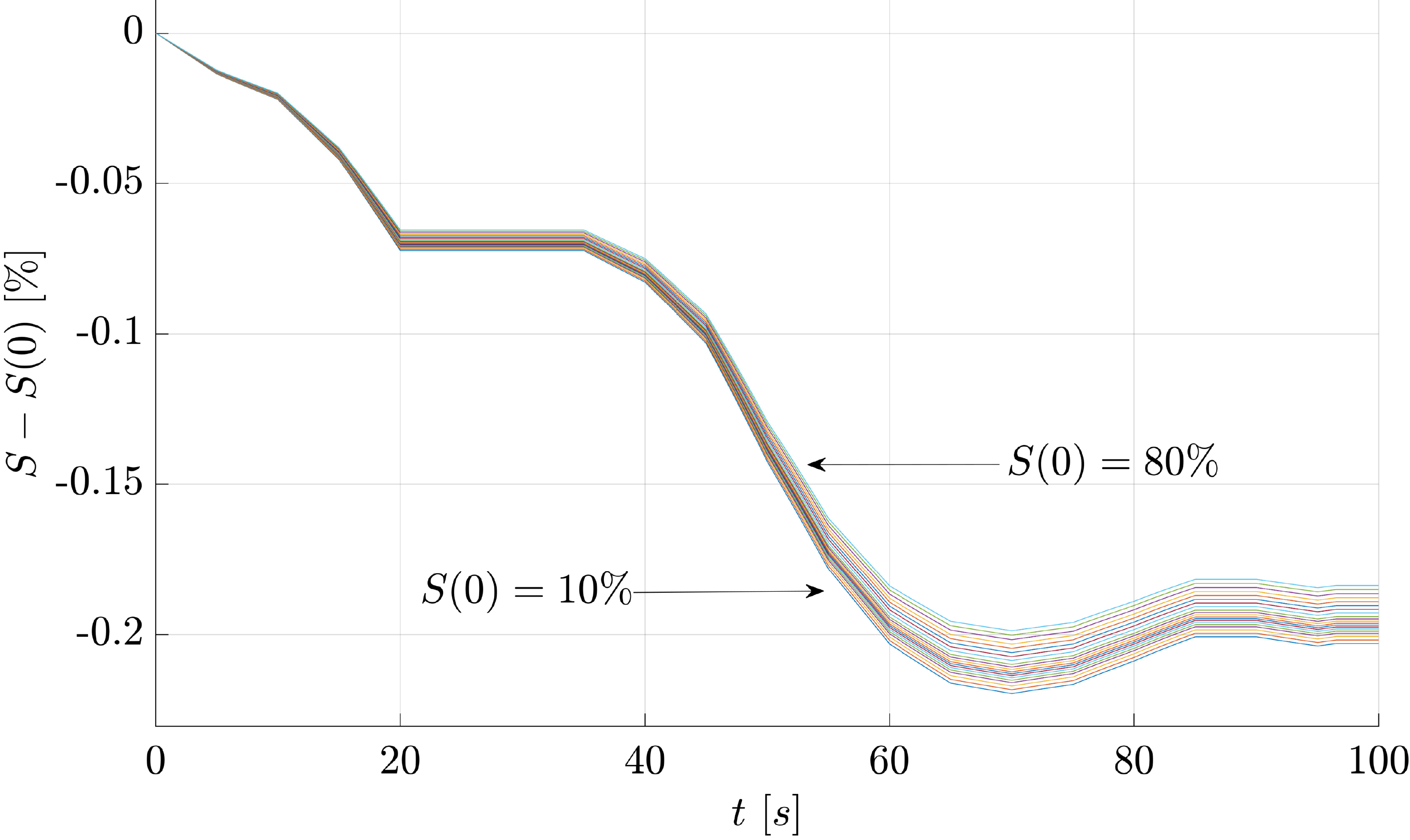}\\(a)} \\
		\parbox[b]{0.24\textwidth}{\centering \includegraphics[width=0.24\textwidth]{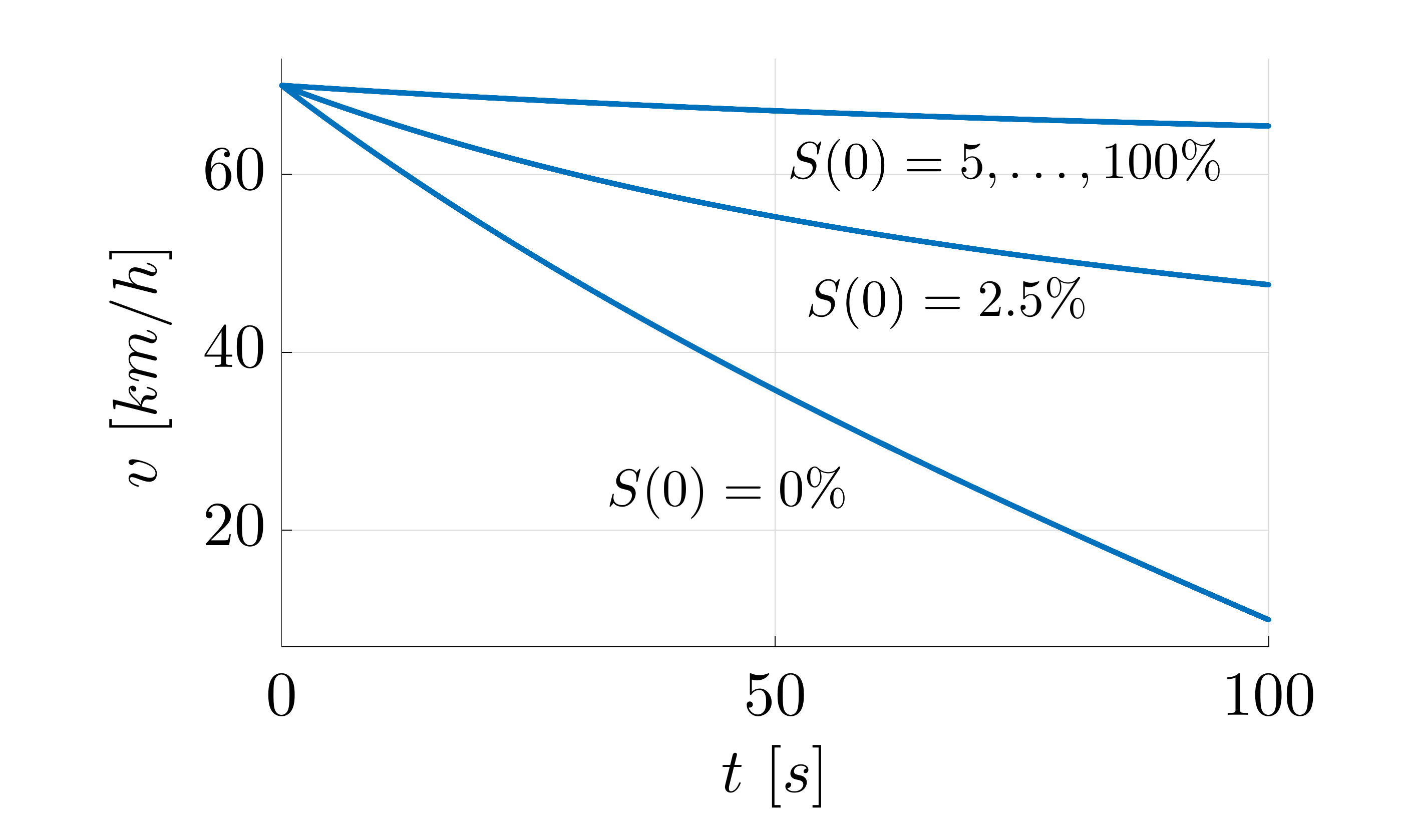}\\(b)}
		\parbox[b]{0.24\textwidth}{\centering \includegraphics[width=0.24\textwidth]{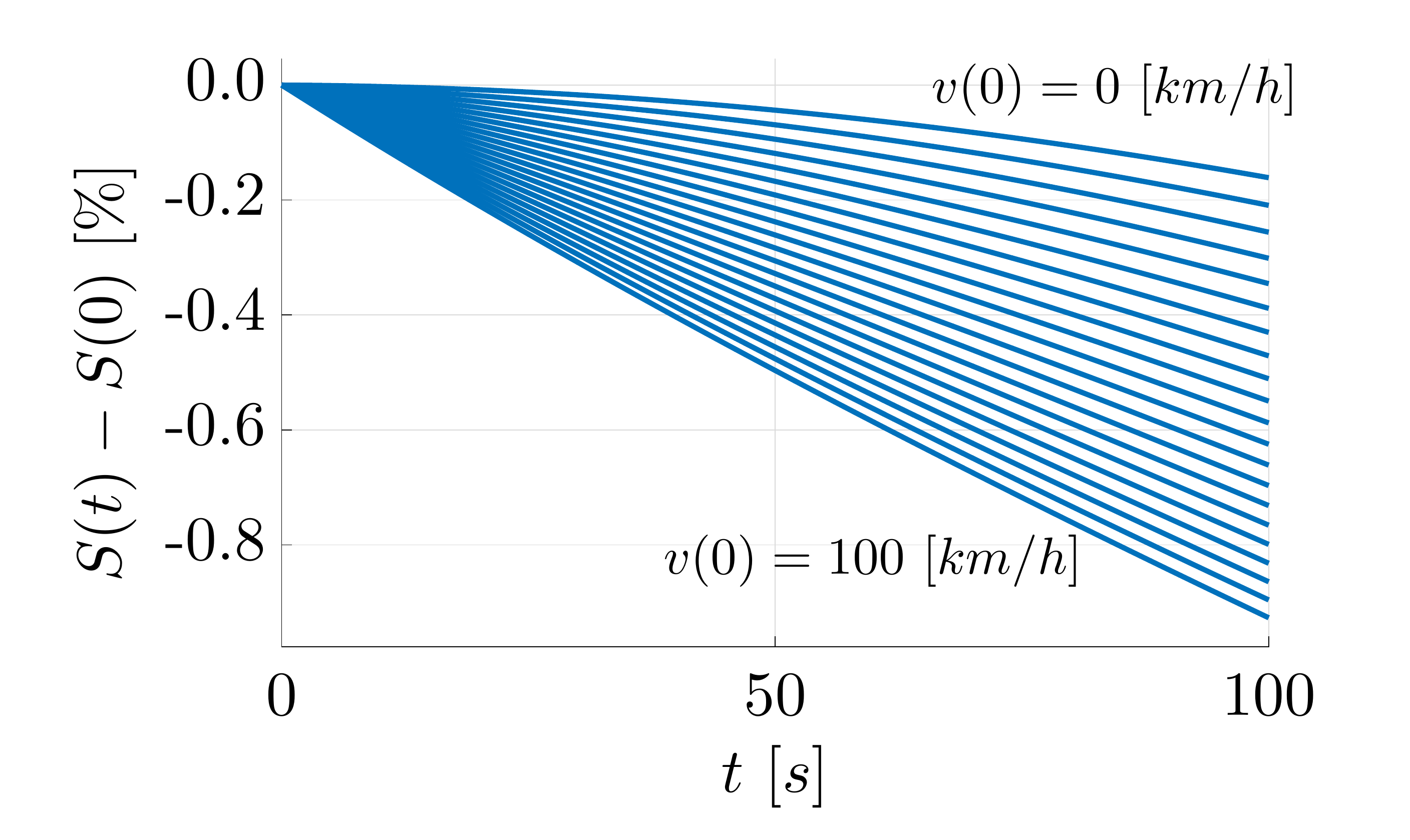}\\(c)}
		\caption{(a) Almost invariance of $S(t)$ with respect to the initial value $S(0)$. (b) Invariance of the velocity $v(t)$ with respect to the initial value $S(0)$ for $S(0) \ge 5\%$. (c) No invariance of the state of charge $S(t)$ with respect to the initial velocity $v(0)$.}
		\label{fig:Invariances}
	\end{figure}
	The more invariances the MOCP possesses (in the sense of Equation~\eqref{eq:Invariance_MOCP}), the fewer problems need to be solved which significantly reduces the computational effort. Hence, we numerically analyse the system in this regard. Since the position $p$ does not occur in the dynamical system \eqref{eq:MOCP_EV_diff}, the dynamics are obviously invariant under translations in $p$. Moreover, when exemplary looking at the velocity $v$ and the state of charge $S$ (cf.~Figures~\ref{fig:Invariances}a and \ref{fig:Invariances}b), we see that, on the one hand, the trajectories are almost invariant for a wide range of translated initial values of the state of charge $S(0)$. 
	Note that this is not a strict invariance. However, as argued in Section \ref{subsec:Algorithm}, we do not require invariances according to Equation~\eqref{eq:Invariance} but according to the weaker condition~\eqref{eq:Invariance_MOCP} which is satisfied much more accurately for the EV application. When looking at Figure~\ref{fig:Invariances}c on the other hand, we observe that the dynamics are clearly not invariant under translations in the initial velocity $v(0)$. After performing the same analysis with regards to the other state variables $U_{d,L}$ and $U_{d,S}$, we can conclude that we only need to define scenarios with respect to the initial velocity $v(0)$ and the active constraints $v_{min}(t)$ and $v_{max}(t)$.
	
	\subsubsection{Constraints}
	\label{subsec:EV_constraints}
	\begin{figure}[h!]
		\centering
		\parbox[b]{0.24\textwidth}{\centering \includegraphics[width=0.24\textwidth]{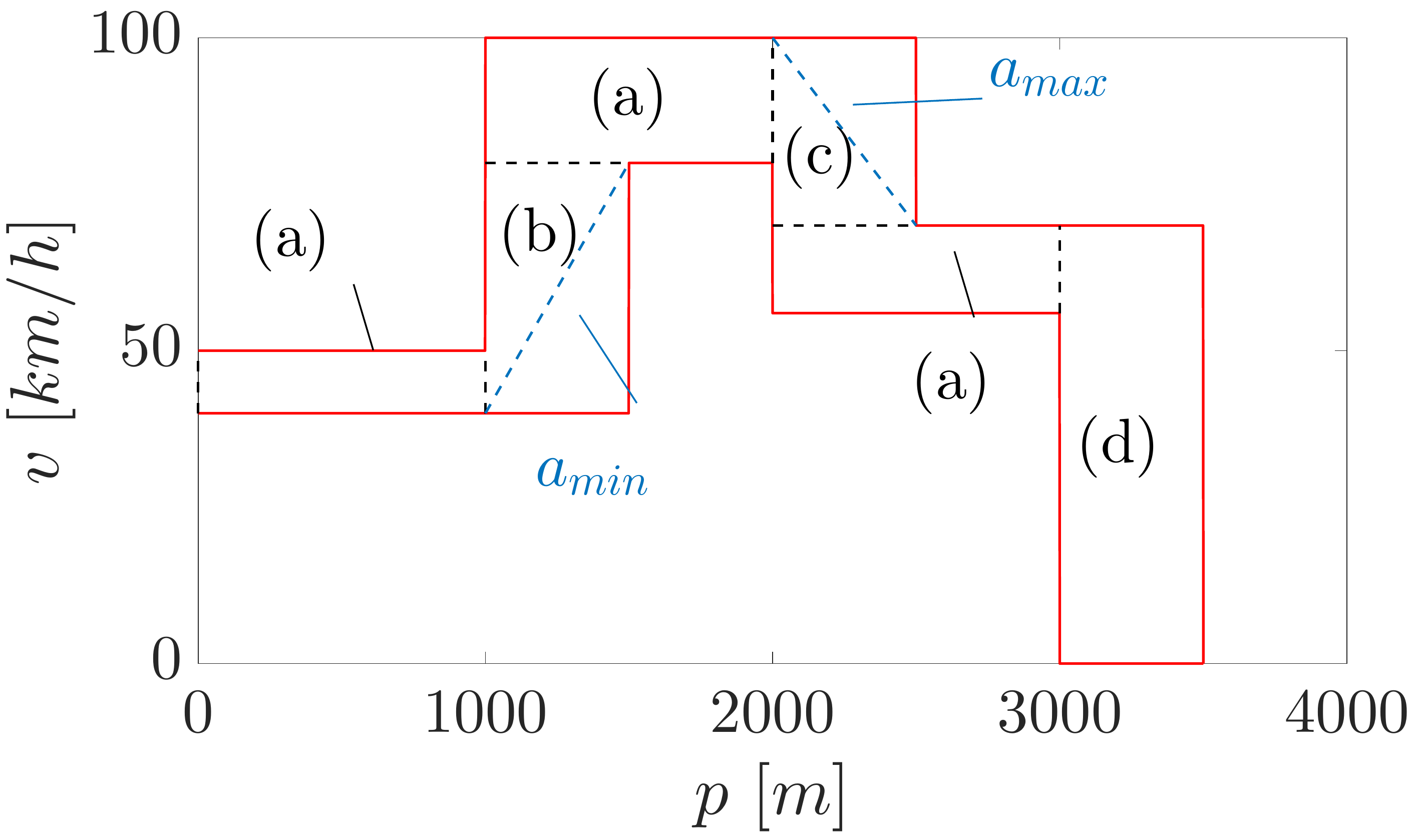}\\(a)}
		\parbox[b]{0.24\textwidth}{\centering \includegraphics[width=0.225\textwidth]{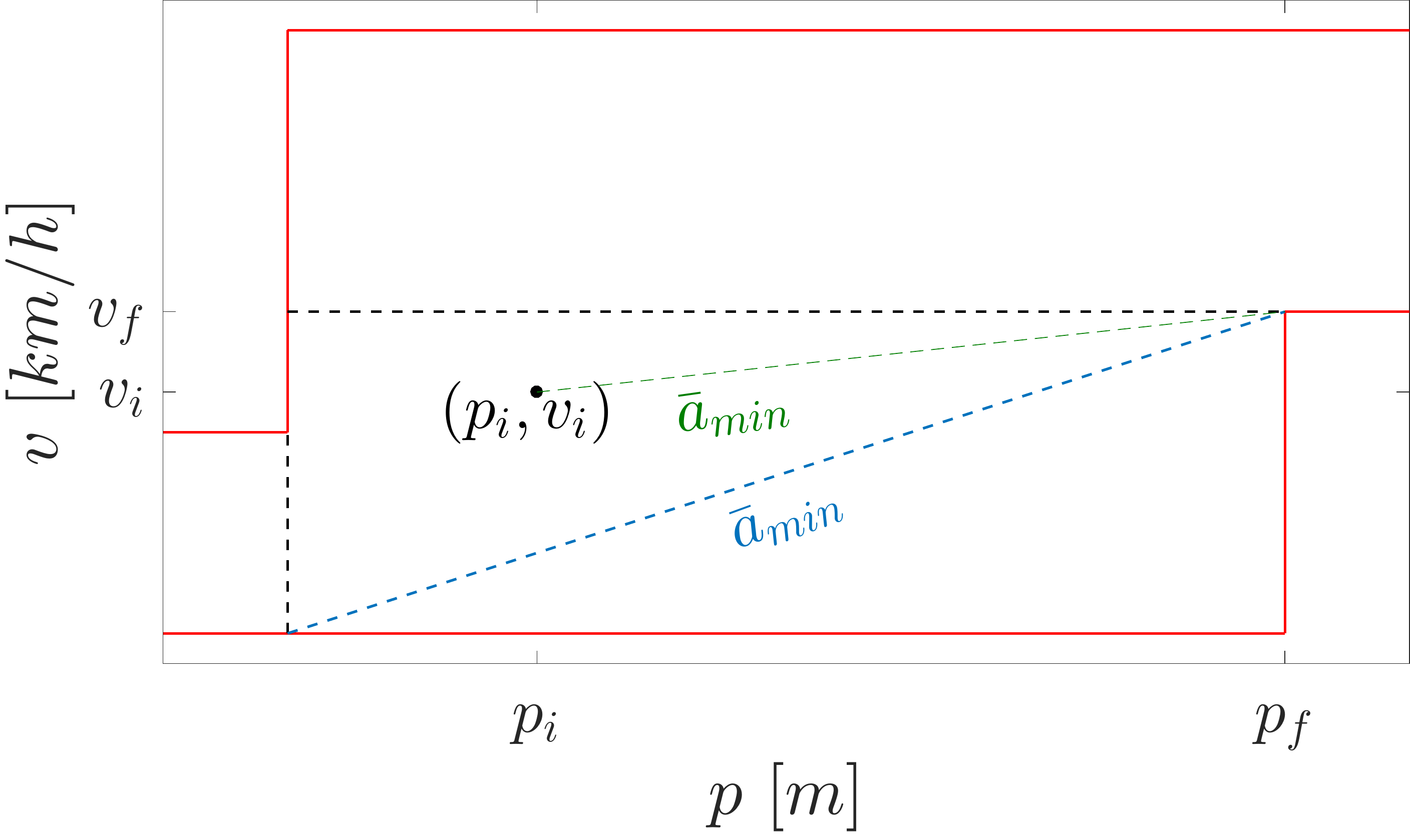}\\(b)}
		\caption{(a) Possible scenarios of boundary conditions. a: constant velocity. b: acceleration. c: deceleration. d: stop sign. (b) Computation of lower bound $\overline{a}_{min}$ for the velocity gradient $dv / dp$.}
		\label{fig:Constraints}
	\end{figure}
\begin{figure}[h!]
		\centering
		\parbox[b]{0.24\textwidth}{\centering \includegraphics[width=0.24\textwidth]{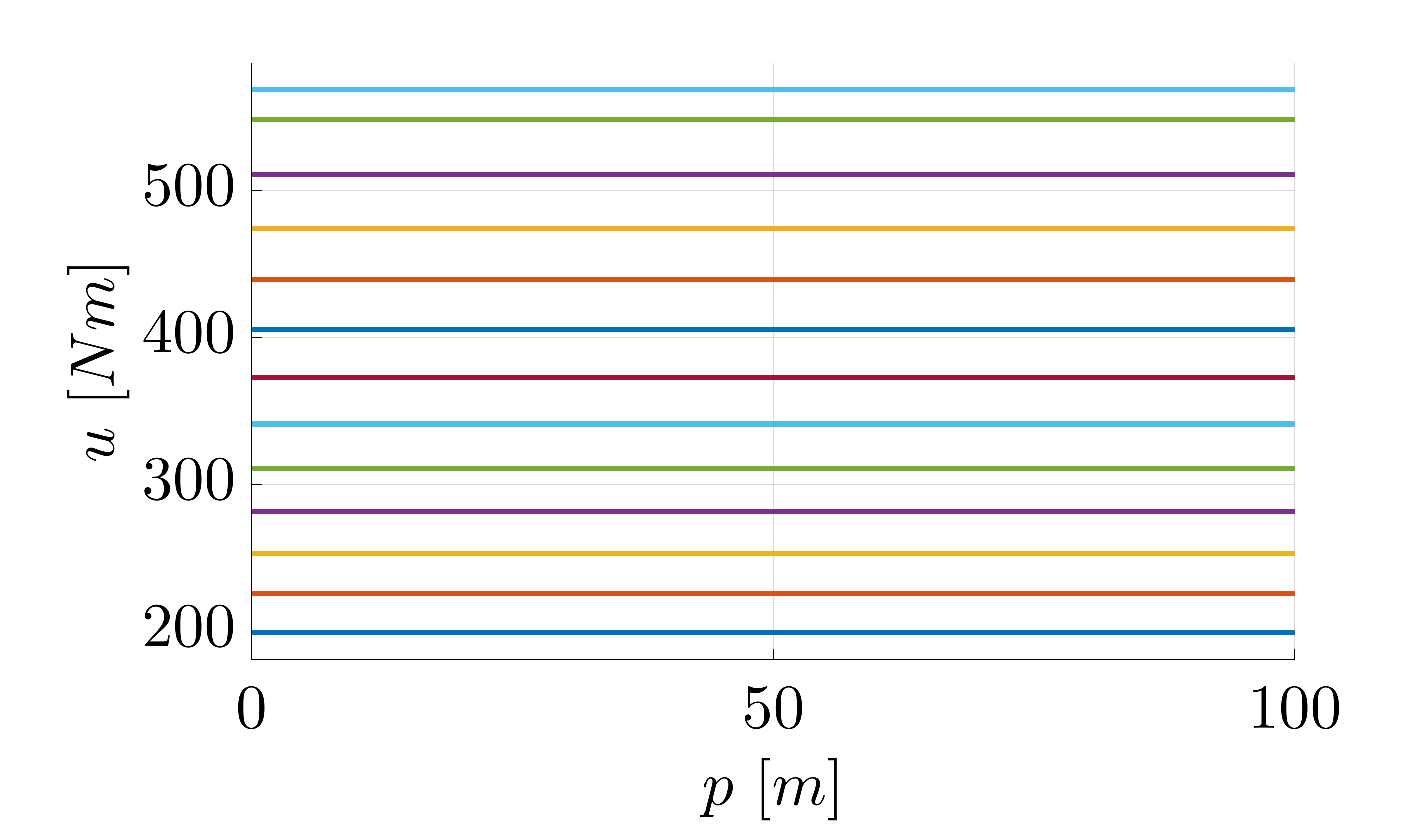}\\(a)}
		\parbox[b]{0.24\textwidth}{\centering \includegraphics[width=0.24\textwidth]{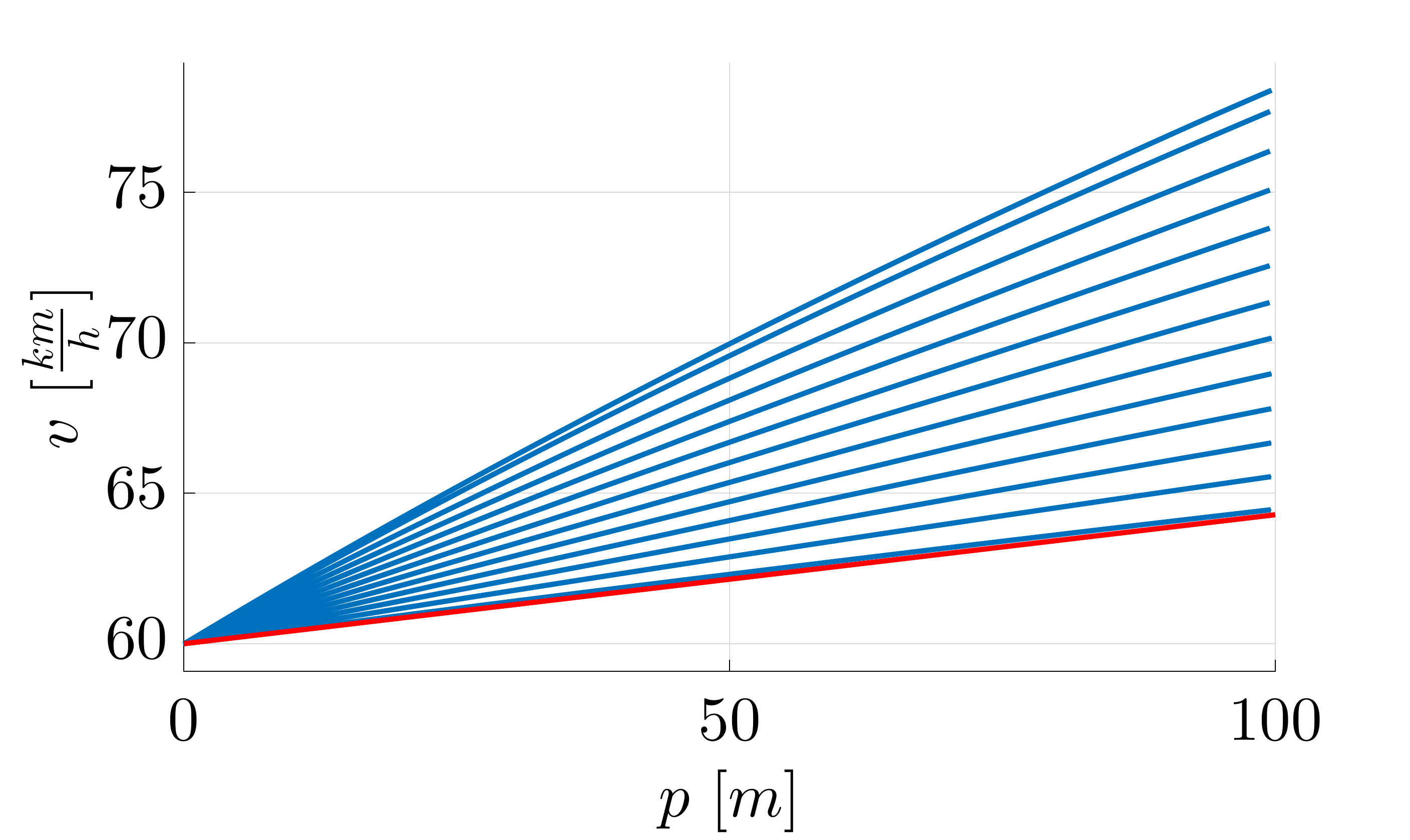}\\(b)}
		\caption{(a) Pareto set for an accelerating scenario with $v(0) = 60\ km/h$ and $\overline{a}_{min} = 0.05\ \frac{km/h}{m}$. (b) The corresponding trajectories of $v(t)$.}
		\label{fig:PS_accel}
	\end{figure}

	A constraint on the velocity is given by the current speed limit $v_{max}(p)$ which depends on the current vehicle position. Since we need to avoid interfering with other vehicles by driving too slow, we define a minimal velocity $v_{min}(p) = 0.8\cdot v_{max}(p)$. (Here we have written the velocities as functions of the position because they are given by the problem formulation this way. In the MOCP, they have to be reformulated as functions of time.) Our \emph{set of feasible states} is now determined by the velocity constraints, i.e.~$v_{min}(t)\le v(t) \le v_{max}(t)$, which determine the different scenarios. We distinguish between four cases (see Figure~\ref{fig:Constraints}a). While the cases constant velocity (box constraints) and stopping ($v=0$ at the stop sign) are easily implemented, we introduce a linear constraint for the scenarios (b) and (c), respectively (see Figure~\ref{fig:Constraints}b) where, depending on the current velocity, a minimal increase $\overline{a}_{min} = (dv/dp)_{min}$ or decrease, respectively, must not be violated. An example is shown in Figure~\ref{fig:PS_accel}, where the Pareto set (\ref{fig:PS_accel}a) and the resulting velocity profiles (\ref{fig:PS_accel}b) are shown for the scenario $v(0) = 60\ km/h$ and $\overline{a}_{min} = 0.05\ \frac{km/h}{m}$. 
	Note that here, we have chosen the control $u$ to be constant over the prediction horizon in order to reduce the numerical effort. As mentioned in Section~\ref{subsec:Algorithm}, we cannot solve an MOCP for every initial condition. Solving an MOCP for every step of $0.1$ in the initial velocity leads to 1727 MOCPs in total.
	
	\subsection{Online Phase: Multiobjective MPC with Paretooptimal Control Primitives}
	\label{subsec:EV_online_phase}
	The online phase is now exactly as described in Section~\ref{subsec:Algorithm}. In each sample time, the current velocity and the active constraints (for the current position) are evaluated in order to determine the valid scenario. The corresponding Pareto set is then selected from the library and according to the weighting parameter $\rho \in [0,1]$ determined by the decision maker, an optimal compromise is chosen which is then applied to the system. On a standard computer, this operation takes in the order of $10^{-3}$ seconds in Matlab.
	
	\section{Results and Discussion}
	\label{sec:Results}
	\begin{figure}[h!]
		\centering
		\includegraphics[width=0.48\textwidth]{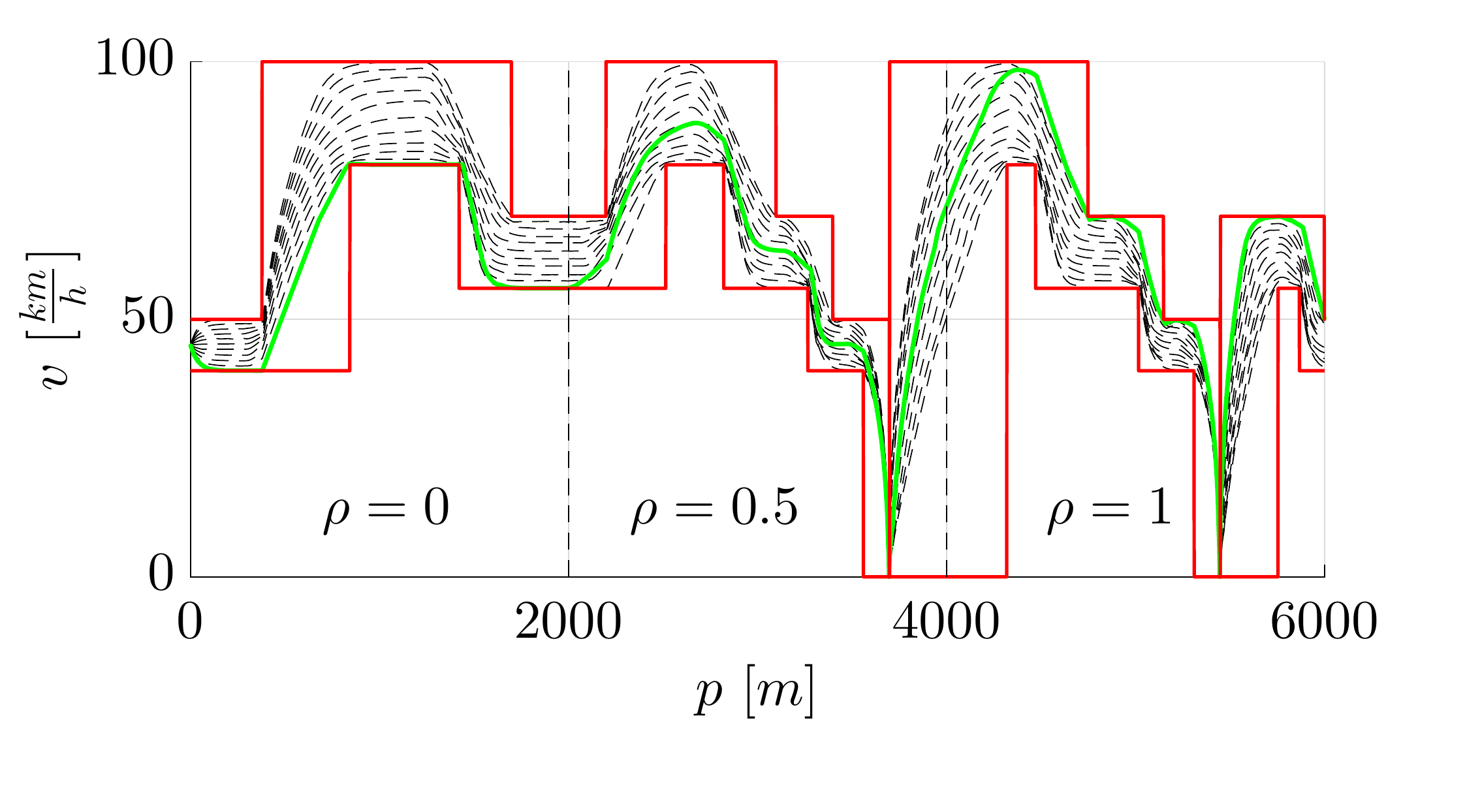}
		\caption{Different trajectories computed by the MPC approach. The dashed lines use a constant weight $\rho$ whereas the green line possesses dynamic weighting ($\rho = 0 / 0.5 / 1.0$, respectively)}
		\label{fig:Road_var_rho}
	\end{figure}
\begin{figure}[h!]
		\centering
		\parbox[b]{0.24\textwidth}{\centering \includegraphics[width=0.24\textwidth]{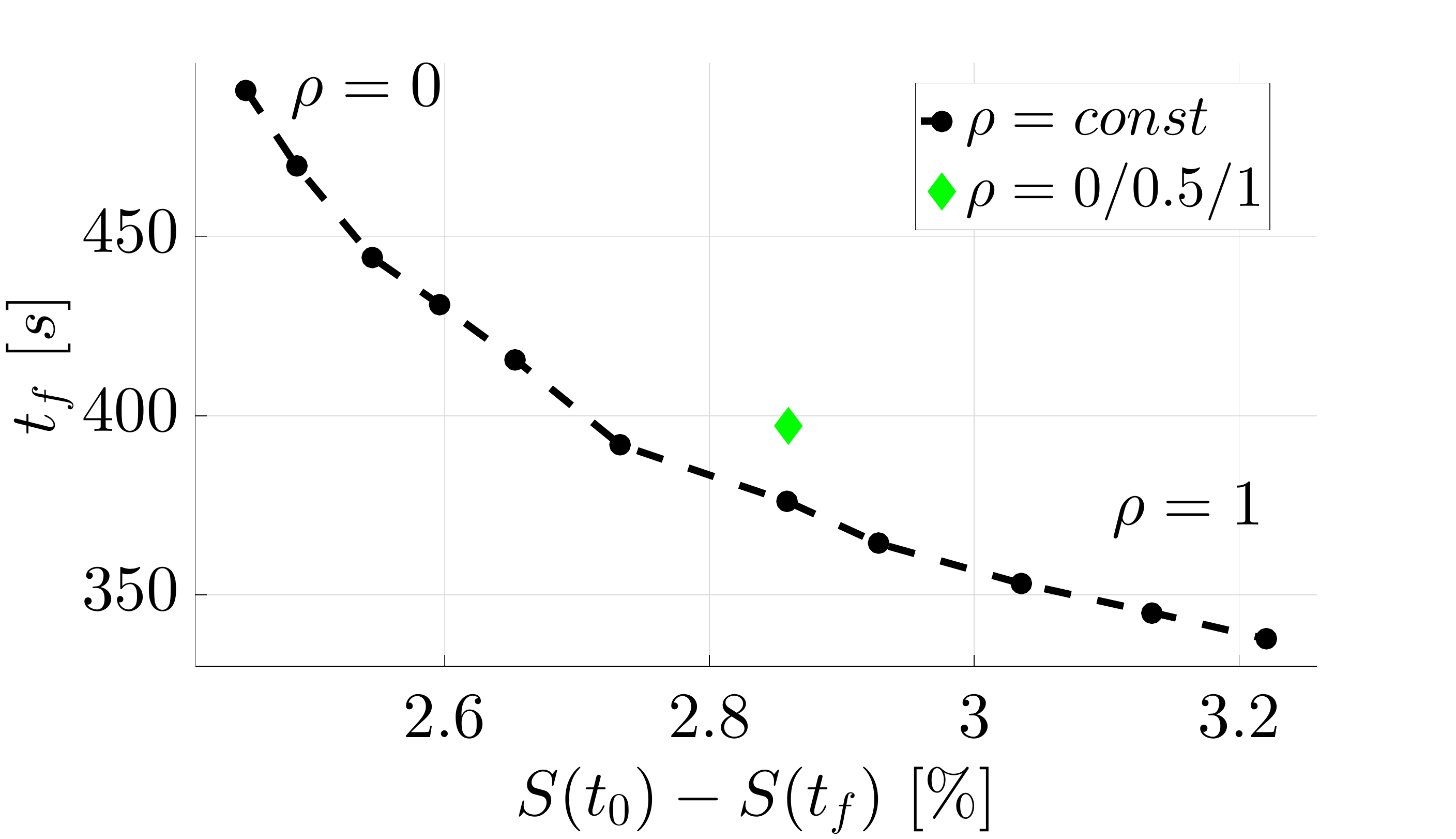}\\(a)}
		\parbox[b]{0.24\textwidth}{\centering \includegraphics[width=0.24\textwidth]{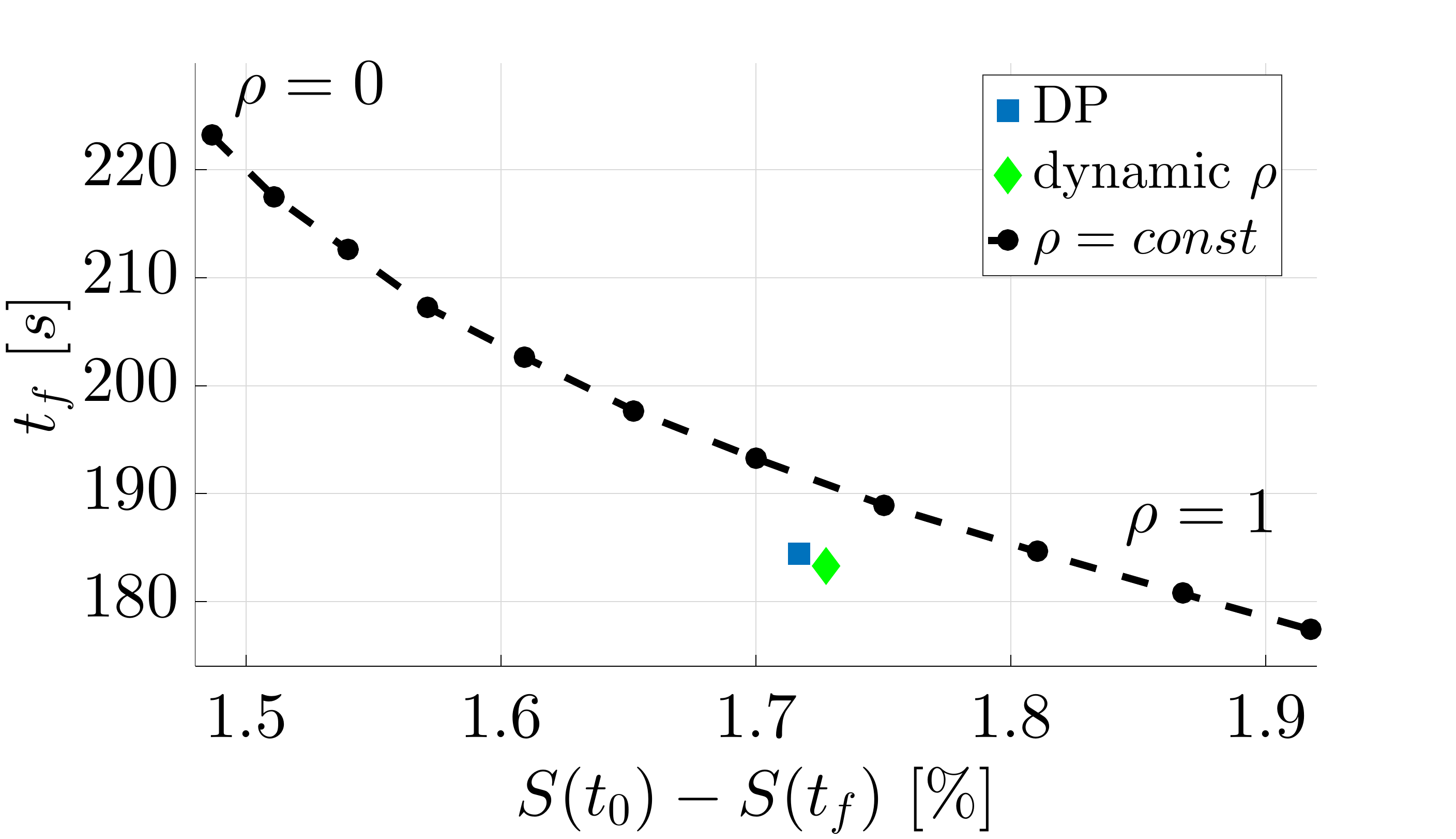}\\(b)}
		\caption{Function values for the scenarios depicted in Figure~\ref{fig:Road_var_rho} and in Figure~\ref{fig:Road_var_rho_dyn} for different weights $\rho$ and in comparison to the Dynamic Programming solution.}
		\label{fig:PF}
	\end{figure}

	In Figure~\ref{fig:Road_var_rho}, several solutions with different weights $\rho$ are shown for an example track including two stop signs. The \emph{set of feasible states} is bounded by the red lines $v_{min}$ and $v_{max}$. The dashed lines correspond to constant weights, varying from $\rho=0$ (energy efficiency) to $\rho = 1$ (high velocity) and the solid green line is a solution where the weighting is changed from 0 over 0.5 to 1 during driving. We clearly see that the vehicle is driving according to the decision maker's preference. This means that we have realised a closed-loop control for which the objectives can be adjusted dynamically. This can either be done manually or by an additional algorithm, which for example takes into account the track, the battery state of charge and the current traffic. The objective function values for the entire track and different values of $\rho$ are depicted in Figure~\ref{fig:PF}a.
	
	In order to evaluate the quality of our solution, we compare it to a control computed via dynamic programming (DP, see~\cite{BD15} for an introduction and \cite{SG09} for the algorithm that is used): For computational reasons, the comparison is performed on a shorter track without stop signs and a relatively coarse discretisation leading to a 100-dimensional problem. In the DP problem, we use a simplified linear model (cf.~\cite{EPS+16}) and the objective is a weighted sum of the MOCP \eqref{eq:MOCP_EV_J}, $J= t_f + \beta E(t_f)$, where $E$ is the consumed energy computed by integrating over the wheel torque and $\beta = 6\cdot10^{-5}$. In Figure~\ref{fig:PF}b, we see that the solution obtained via DP is superior to our MPC approach. This is not surprising since in MPC, we only consider a finite horizon such that the results are at best suboptimal (\cite{GP11}), whereas the entire track is considered at once in DP. Consequently, the DP algorithm is not real-time applicable and does not possess feedback behaviour. Additionally, we have until now only considered constant torques over the prediction horizon in our approach. We intend to refine the discretisation in future work and expect an improved performance. 
	
	\begin{figure}[h!]
		\centering
		\includegraphics[width=0.48\textwidth]{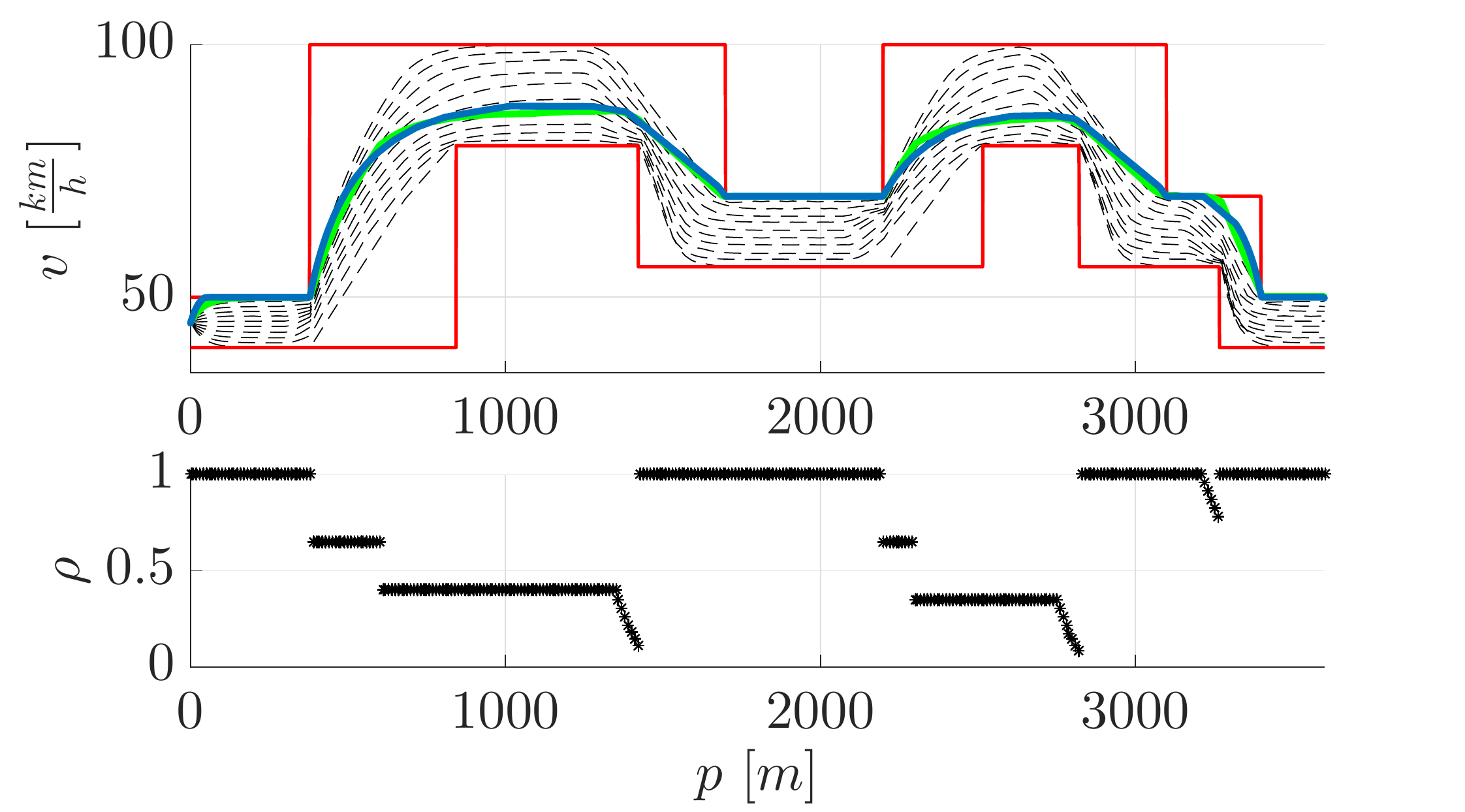}
		\caption{Validation of the approach versus a Dynamic Programming solution (blue). Green line: dynamic weighting according to the lower plot.}
		\label{fig:Road_var_rho_dyn}
	\end{figure}
	When using a simple, manually tuned heuristic for the preference $\rho$ instead of fixed values (larger values for $\rho$ at low velocities, lower values at high velocities and linear changes in $\rho$ when approaching braking manoeuvres, see Figure~\ref{fig:Road_var_rho_dyn}, bottom), we see that we can improve the quality of our solution significantly which is now comparable to the global optimum obtained by DP. We see in Figures~\ref{fig:Road_var_rho_dyn} (top) and \ref{fig:PF}b, respectively, that the resulting trajectories as well as the function values $J_1$ and $J_2$ almost coincide. By this, we obtain two different ways to utilise the results. On the one hand, a decision maker can select the preference according to his wishes and on the other hand, $\rho$ can be determined by a heuristic, leading to solutions of a quality comparable to the global optimum.
	
	\section{Conclusion}
	\label{sec:Conclusion}
	We present an algorithm for MPC of non-linear dynamical systems with respect to multiple criteria. The algorithm utilises elements from economic and explicit MPC, multiobjective optimal control and motion planning. According to a decision maker's preference, the system is controlled in real-time with respect to an optimal compromise between conflicting objectives. Using a simple heuristic for the weighting factor $\rho$, we obtain solutions of equivalent quality compared to a global optimum computed by open loop DP. In the future, we intend to analyse the proposed method from a more theoretical point of view, addressing questions concerning feasibility and stability for systems where these aspects are critical. Furthermore, we want to improve our control strategies by developing intelligent heuristics for the preference weighting function $\rho$.
	\\~
	
	\textbf{Acknowledgement:} This research was funded by the German Federal Ministry of Education and Research (BMBF) within the Leading-Edge Cluster \emph{Intelligent Technical Systems OstWestfalenLippe (it's OWL)}.
	
	\bibliographystyle{unsrt}
	\setlength{\bibsep}{0pt plus 0.3ex}
	\bibliography{arXiv_bibliography}

\begin{thebibliography}{10}

\bibitem{Ehr05}
M.~Ehrgott.
\newblock {\em {Multicriteria optimization}}.
\newblock Springer Berlin Heidelberg New York, 2005.

\bibitem{CLV07}
C.~A. {Coello Coello}, G.~B. Lamont, and D.~A. van Veldhuizen.
\newblock {\em {Evolutionary Algorithms for Solving Multi-Objective Problems}},
  volume~2.
\newblock Springer New York, 2007.

\bibitem{SWO+13}
O.~Sch{\"u}tze, K.~Witting, S.~Ober-Bl{\"o}baum, and M.~Dellnitz.
\newblock {Set Oriented Methods for the Numerical Treatment of Multiobjective
  Optimization Problems}.
\newblock In Emilia Tantar~et al., editor, {\em EVOLVE - A Bridge between
  Probability, Set Oriented Numerics and Evolutionary Computation}, volume 447
  of {\em Studies in Computational Intelligence}, pages 187--219. Springer
  Berlin Heidelberg, 2013.

\bibitem{Mac02}
J.~M. Maciejowski.
\newblock {\em {Predictive Control: With Constraints}}.
\newblock Prentice Hall, Harlow, England, 2002.

\bibitem{GP11}
L.~Gr{\"u}ne and J.~Pannek.
\newblock {\em {Nonlinear model predictive control}}.
\newblock Springer, 2011.

\bibitem{QB97}
S.~J. Qin and T.~A. Badgwell.
\newblock {An overview of industrial model predictive control technology}.
\newblock In {\em AIChE Symposium Series}, volume~93, pages 232--256. American
  Institute of Chemical Engineers, 1997.

\bibitem{LC97}
J.~H. Lee and B.~Cooley.
\newblock {Recent advances in model predictive control and other related
  areas}.
\newblock In {\em AIChE Symposium Series}, volume~93, pages 201--216. American
  Institute of Chemical Engineers, 1997.

\bibitem{EPS+16}
J.~Eckstein, K.~Sch{\"a}fer, S.~Peitz, P.~Friedel, S.~Ober-Bl{\"o}baum, and
  M.~Dellnitz.
\newblock {A Comparison of two Predictive Approaches to Control the
  Longitudinal Dynamics of Electric Vehicles}.
\newblock {\em Procedia Technology}, 26:465--472, 2016.

\bibitem{RA09}
J.~B. Rawlings and R.~Amrit.
\newblock Optimizing process economic performance using model predictive
  control.
\newblock In {\em Nonlinear model predictive control}, pages 119--138.
  Springer, 2009.

\bibitem{DAR11}
M.~Diehl, R.~Amrit, and J.~B. Rawlings.
\newblock A lyapunov function for economic optimizing model predictive control.
\newblock {\em IEEE Transactions on Automatic Control}, 56(3):703--707, 2011.

\bibitem{AB09}
A.~Alessio and A.~Bemporad.
\newblock {A Survey on Explicit Model Predictive Control}.
\newblock In Lalo Magni, Davide~Martino Raimondo, and Frank Allg{\"o}wer,
  editors, {\em Nonlinear Model Predictive Control: Towards New Challenging
  Applications}, pages 345--369. Springer Berlin Heidelberg, 2009.

\bibitem{FDF05}
E.~Frazzoli, M.~A. Dahleh, and E.~Feron.
\newblock {Maneuver-Based Motion Planning for Nonlinear Systems with
  Symmetries}.
\newblock {\em IEEE Transactions on Robotics}, 21(6):1077--1091, 2005.

\bibitem{Kob08}
M.~Kobilarov.
\newblock {\em {Discrete geometric motion control of autonomous vehicles}}.
\newblock PhD thesis, University of Southern California, 2008.

\bibitem{FOK12}
K.~Fla{\ss}kamp, S.~Ober-Bl{\"o}baum, and M.~Kobilarov.
\newblock {Solving Optimal Control Problems by Exploiting Inherent Dynamical
  Systems Structures}.
\newblock {\em Journal of Nonlinear Science}, 22(4):599--629, 2012.

\bibitem{BP09}
A.~Bemporad and D.~{Mu{\~{n}}oz de la Pe{\~{n}}a}.
\newblock {Multiobjective model predictive control}.
\newblock {\em Automatica}, 45(12):2823--2830, 2009.

\bibitem{ZFT12}
V.~M. Zavala and A.~Flores-Tlacuahuac.
\newblock {Stability of multiobjective predictive control: A utopia-tracking
  approach}.
\newblock {\em Automatica}, 48(10):2627--2632, 2012.

\bibitem{DEF+16}
M.~Dellnitz, J.~Eckstein, K.~Fla{\ss}kamp, P.~Friedel, C.~Horenkamp,
  U.~K{\"o}hler, S.~Ober-Bl{\"o}baum, S.~Peitz, and S.~Tiemeyer.
\newblock {Multiobjective Optimal Control Methods for the Development of an
  Intelligent Cruise Control}.
\newblock In G.~Russo~et al., editor, {\em Progress in Industrial Mathematics
  at ECMI 2014 (to appear)}, 2016.

\bibitem{NW06}
J.~Nocedal and S.~J. Wright.
\newblock {\em {Numerical Optimization}}.
\newblock Springer Science {\&} Business Media, 2006.

\bibitem{DEF+14}
M.~Dellnitz, J.~Eckstein, K.~Fla{\ss}kamp, P.~Friedel, C.~Horenkamp,
  U.~K{\"o}hler, S.~Ober-Bl{\"o}baum, S.~Peitz, and S.~Tiemeyer.
\newblock {Development of an Intelligent Cruise Control Using Optimal Control
  Methods}.
\newblock In {\em Procedia Technology}, volume~15, pages 285--294. Elsevier,
  2014.

\bibitem{BD15}
R.~E. Bellmann and S.~E. Dreyfus.
\newblock {\em {Applied dynamic programming}}.
\newblock Princeton University Press, 2015.

\bibitem{SG09}
O.~Sundstr{\"o}m and L.~Guzzella.
\newblock {A generic dynamic programming Matlab function}.
\newblock In {\em 2009 IEEE Control Applications, (CCA) \& Intelligent
  Control,(ISIC)}, pages 1625--1630, 2009.

\end{thebibliography}
\end{document}